\documentclass{amsart}

\usepackage[colorlinks,
linkcolor=blue,
anchorcolor=blue,
citecolor=blue]{hyperref}

\usepackage{verbatim}
\newtheorem{theorem}{Theorem}[section]
\newtheorem{lemma}{Lemma}[section]
\newtheorem{proposition}{Proposition}[section]
\newtheorem{corollary}{Corollary}[section]
\newtheorem{question}{Question}
\newtheorem{conjecture}[theorem]{Conjecture}
\newtheorem{problem}{Problem}

\theoremstyle{definition}
\newtheorem{definition}{Definition}[section]
\newtheorem{example}{Example}[section]

\newtheorem{remark}{Remark}[section]

\numberwithin{equation}{section}

\DeclareMathOperator{\spa}{Span}
\DeclareMathOperator{\ind}{ind}
\DeclareMathOperator{\tr}{Tr}

\begin{document}
\title{Hilbert-Schmidtness of the $M_{\theta,\varphi}$-type submodules }
\author{Chao Zu}
\address{School of Mathematical Sciences, Dalian University of Technology,
Dalian, Liaoning, 116024, P. R. China}
\email{zuchao@dlut.edu.cn}

\author{Yufeng Lu}
\address{School of Mathematical Sciences, Dalian University of Technology,
Dalian, Liaoning, 116024, P. R. China}
\email{lyfdlut@dlut.edu.cn}


\begin{abstract}
Let $\theta(z),\varphi(w)$ be two nonconstant inner functions and $M$ be a submodule in $H^2(\mathbb{D}^2)$. Let $C_{\theta,\varphi}$ denote the composition operator on $H^2(\mathbb{D}^2)$ defined by $C_{\theta,\varphi}f(z,w)=f(\theta(z),\varphi(w))$, and $M_{\theta,\varphi}$ denote the submodule $[C_{\theta,\varphi}M]$, that is, the smallest submodule containing $C_{\theta,\varphi}M$. Let $K^M_{\lambda,\mu}(z,w)$ and $K^{M_{\theta,\varphi}}_{\lambda,\mu}(z,w)$ be the reproducing kernel of $M$ and  $M_{\theta,\varphi}$, respectively. By making full use of the positivity of certain de Branges-Rovnyak kernels, we prove that
\[K^{M_{\theta,\varphi}}= K^M \circ B~ \cdot R,\]
where $B=(\theta,\varphi)$, $R_{\lambda,\mu}(z,w)=\frac{1-\overline{\theta(\lambda)}\theta(z)}{1-\bar{\lambda}z} \frac{1-\overline{\varphi(\mu)}\varphi(w)}{1-\bar{\mu}w}$. This implies that $M_{\theta,\varphi}$ is a Hilbert-Schmidt submodule if and only if $M$ is. Moreover, as an application, we prove that the Hilbert-Schmidt norms of submodules $[\theta(z)-\varphi(w)]$ are uniformly bounded.
\end{abstract}

\maketitle

\section{Introduction}
The Hardy space over the bidisk $H^2(\mathbb{D}^2)$, under action defined by multiplication of functions, is a module over the polynomial ring $\mathbb{C}[z,w]$. Submodules
of $H^2(\mathbb{D}^2)$ have complicated yet intriguing structure, which has attracted  continuing efforts in search of an elucidation. Since the 60's, a great amount of work has been done. One prevalent idea in this area is to define equivalence relations among submodules and study equivalence classes. Various kinds of questions regarding unitary equivalence thus arise \cite{CG}. A most notable fact is the rigidity phenomenon discovered by R. Douglas, V. Paulsen,  C.-H. Sah and K. Yan in \cite{DPSY}. Another important line of research is to consider the so-call Hilbert-Schmidt submodules. Such submodules, on the one hand, are broad enough to include almost all known examples, and on the other hand, they are good enough so that some fine analysis tools apply. The concept of ``Hilbert-Schmidt submodule" initially originated from R. Yang's generalization of the classical Berger-Shaw theorem  \cite{Ya1999}, which says the self-commutator of a multicyclic hyponormal operator is in the trace class. Here, an operator $A$ on a Hilbert space $H$ is said to be $n$-multicyclic if there exist $n$ vectors $g_1,\cdots,g_n$ in $H$ for which the linear manifold spanned by the vectors $r(A)g_i$, where $r$ runs through the rational functions with poles off the spectra set $\sigma(A)$, is dense in $H$. An operator $A$ is hyponormal if its self-commutator $[A^*,A]=A^*A-AA^*$ is nonnegative.

\begin{theorem}[Berger-Shaw Theorem, \cite{Berger}]
If $T$ is a $n$-multicyclic hyponormal operator, then its self-commutator $[T^*,T]$ is of trace class and
\[  \tr[T^*,T]\leq \frac{n}{\pi}~ \mu(\sigma(T)) , \]
where $\mu$ is the Lebesque measure on the plane.
\end{theorem}
It is interesting to find a multivariate analogue of this theorem. This naturally leads to the question: are there natural definitions of subnormal pairs and self-commutator which admit a generalization of the Berger-Shaw theorem? Let $R_z, R_w$ be the multiplication of the coordinate functions $z$ and $w$ on submodule $M$, and let
$S_z, S_w$ be the actions of $z$ and $w$ on quotient module $H^2(\mathbb{D}^2)\ominus M$. Inspired by the work of R. Curto, P. Muhly and K. Yan \cite{Curto} and
the ideas of R. Douglas \cite{Douglas1,Douglas2}, R. Yang studied the Hilbert-Schmidtness of the product $[R_z^*,R_z][R_w^*,R_w]$ and the cross commutators $[R_z^*,R_w]$, $[S_z^*,S_w]$, and defined the \textbf{numerical invariants} for submodules:
\[ \Sigma_0^M:=\|[R_z^*,R_z][R_w^*,R_w]\|_{H.S.}^2,~~\Sigma_1^M:=\|[R_z^*,R_w]\|_{H.S.}^2.\]
A submodule is said to be \textbf{Hilbert-Schmidt} if $\Sigma_0^M$ and $\Sigma_1^M$ are finite. It was shown by R.Yang \cite{Ya3} that $M$ is Hilbert-Schmidt if $\sigma_e(M):=\sigma_e(S_z)\cap \sigma_e(S_w) \neq \overline{\mathbb{D}}$ (which is a very mild condition!). Moreover, R. Yang proved that every finitely generated polynomial submodule is  Hilbert-Schmidt and made the following conjecture.
\begin{conjecture}[R. Yang, 1999, \cite{Ya1999}]
Every finitely generated submodule in $H^2(\mathbb{D}^2)$ is Hilbert-Schmidt.
\end{conjecture}
In \cite{Luo}, S. Luo, K. Izuchi and R. Yang confirmed the conjecture for submodules containing function $z-\varphi(w)$, where $\varphi(w)$ is any finite Blaschke product. Recently, we extended this result to the submodules containing function $\theta(z)-\varphi(w)$, where $\theta,\varphi$ are two finite Blaschke products \cite{LYZ}. This also inspires us to study the connection between the submodule $M_{\theta,\varphi}$ and the submodule $M$.

Further, R. Yang also proposed the following more general conjecture.
\begin{conjecture}\label{rank-one}
Let $\mathcal{S}$ be the set of all rank-one submodules in $H^2(\mathbb{D}^2)$, then
\[ \sup_{M\in \mathcal{S}} ~\Sigma_0^M+\Sigma_1^M < \infty. \]
\end{conjecture}
In particular, as polynomial submodules are always Hilbert-Schmidt, we can first consider above conjecture for polynomial submodules:
\begin{question}
Does it hold that
\[ \sup_{q\in \mathbb{C}[z,w]} \Sigma_0^{[q]}+\Sigma_1^{[q]}<\infty? \]
\end{question}

This paper focuses primarily on  the Hilbert-Schmidtness of submodules $M_{\theta,\varphi}$ and their norm estimation. In section 3, by making full use of the positivity of certain kernel functions (inspired by the work of M. Jury \cite{Jury} and the work of C. Chu \cite{Chu}), we prove that for any submodule $M$, the reproducing kernel of $M_{\theta,\varphi}$ has an intriguing relationship with the reproducing kernel of $M$:
\[K^{M_{\theta,\varphi}}= K^M \circ B~ \cdot R,\]
where $B=(\theta,\varphi)$, $R_{\lambda,\mu}(z,w)=\frac{1-\overline{\theta(\lambda)}\theta(z)}{1-\bar{\lambda}z} \frac{1-\overline{\varphi(\mu)}\varphi(w)}{1-\bar{\mu}w}$ (Theorem \ref{thm3.1}). Then combined with the Littlewood's Subordination Principle, we show that  the submodule $M_{\theta,\varphi}$ is Hilbert-Schmidt if and only if the submodule $M$ is.
In section 4, we extend the numerical invariants $\Sigma_0^M, \Sigma_1^M$ to the numerical invariant functions $\Sigma_0^M(a,b), \Sigma_1^M(a,b)$ and study their properties. We get some useful calculation formulas for $\Sigma_i^M(a,b)$ and find an interesting connection between the numerical invariant functions of $M_{\theta,\varphi}$ and that of $M$  (Theorem \ref{thm4.6}, Corollary \ref{thm4.8}):
\[ \Sigma_0^{M_{\theta,\varphi}}(a,b)=\Sigma_0^M(\theta(a),\varphi(b)),~~~ \Sigma_1^{M_{\theta,\varphi}}(a,b)=\Sigma_1^M(\theta(a),\varphi(b)).\]
Moreover, as an application, in section 5, we prove that the Hilbert-Schmidt norms of submodules $[\theta(z)-\varphi(w)]$ are uniformly bounded(Theorem \ref{thm5.3}), which partially confirms Conjecture \ref{rank-one} .
\section{Preliminaries}
\subsection{Background and notations}
Let $\mathbb{D} = \{\lambda\in \mathbb{C} : ~|\lambda| < 1 \}$ be the unit disk in the complex plane $\mathbb{C}$ and let $\mathbb{T}$ be the unit circle. Let $\mathbb{C}[z,w]$  denote the polynomial ring over $\mathbb{C}$ with variables $z$ and $w$ and $L^2 = L^2(\mathbb{T}^2)$ denote the Hilbert space of all square integrable functions with respect to the normalized Lebesgue measure $d\sigma$ on $\mathbb{T}^2$. $L^\infty=L^\infty(\mathbb{T}^2)$ denotes the commutative Banach algebra of all
essentially bounded measurable functions with respect to the measure $d\sigma$ on $\mathbb{T}^2$. $H^2 = H^2(\mathbb{D}^2)$ will denote the Hardy space over the bidisk $\mathbb{D}^2$. It consists of all analytic functions on $\mathbb{D}^2$ satisfying the following condition:
\[\|f\|_2^2:=\sup_{0\leq r<1} \int_{\mathbb{T}}\int_{\mathbb{T}} |f(rz,rw)|^2 d\sigma <\infty.\]
 Let $H^2_z$ (resp. $H^2_w$)
denote the usual one-variable Hardy space over $\mathbb{D}$ of the variable $z$ (resp. $w$) with
respect to the normalized Lebesgue measure $|dz|$ (resp. $|dw|$) on $\mathbb{T}$. Then it is well
known that $H^2 =H^2_z\otimes H_w^2 $. $H^\infty = H^\infty(\mathbb{D}^2)$ will denote the commutative Banach algebra consisting
of all bounded analytic functions on $\mathbb{D}^2$ with the norm $\|f\|_{\infty}=\sup_{(\lambda,\mu)\in \mathbb{D}^2} |f(\lambda,\mu)|$.
The Toeplitz operator $T_\phi$ on $H^2$, corresponding to a bounded analytic function $\phi$, acts as the multiplication operator by $\phi$ on $H^2$. A closed subspace $M$ of $H^2(\mathbb{D}^2)$ is called a submodule if $M$ is invariant under $T_z$ and $T_w$. If $M$ is a submodule of $H^2(\mathbb{D}^2)$ and $N=H^2(\mathbb{D}^2)\ominus M$, then $N$ is invariant under $T_z^*$ and $T_w^*$, and it is called a quotient module of $H^2(\mathbb{D}^2)$. For any subset $X\subset H^2(\mathbb{D}^2)$, let
\[[X]:=\overline{\spa \{~\mathbb{C}[z,w]\cdot X~\}}^{\|\cdot\|_2}  \]
denote the submodule generated by $X$. For example, $[h]$ is the submodule generated by function $h$. 
A subset $X$ of submodule $M$ is said to be a generating set of $M$ if $[X]=M$. The minimal cardinality of such generating sets for a submodule $M$ is called the rank of $M$.

For two inner functions $\theta(z),\varphi(w)$, let $C_{\theta,\varphi}$ be the composition operator defined by
\[C_{\theta,\varphi}f(z,w)=f(\theta(z),\varphi(w)),~~~~~~~f\in ~Hol(\mathbb{D}^2).\]

If $\varphi$ an analytic mapping on the disk and $w\neq\varphi(0)$ is a point of the plane, let $z_j$ be the points in the disk for which $\varphi(z)=w$, counting multiplicities. The \textbf{Nevanlinna counting function} is defined as
\[N_\varphi(w)=\sum_{j} \log \frac{1}{|z_j|},\]
where we understand that $N_\varphi(w)=0$ for $w$ not in $\varphi(\mathbb{D})$. If $\varphi$ is any analytic map of the disk into itself,  we always have the \textbf{Littlewood's Inequality} (Theorem 2.29, \cite{Cowen})
\[ N_\varphi(w)\leq \log \Big|\frac{1-\overline{\varphi(0)}w}{\varphi(0)-w} \Big| \]  
for $w$ in $\mathbb{D}\backslash \{\varphi(0)\}$. When $\varphi$ is inner, the inequality in Littlewood's Inequality becomes equality, almost everywhere.
\begin{lemma}\label{Nevanlina}
If $\varphi$ is an inner function then
\[N_\varphi(w)= \log \Big|\frac{1-\overline{\varphi(0)}w}{\varphi(0)-w} \Big| \]
for all $w$ in $\mathbb{D}$ outside a set of area measure $0$.
\end{lemma}
The connection between composition operators and counting functions arises directly from the following change of variable formula due to J. H. Shapiro\cite{Shapiro}.
\begin{theorem}\label{Shapiro}
If $f$ is analytic in the unit disk and  $\varphi$ is a non-constant analytic mapping of $\mathbb{D}$ into itself, then
 \[ \|f\circ\varphi\|_2^2-|f(\varphi(0))|^2=2 \int_{\mathbb{D}}|f'(w)|N_\varphi(w) \frac{dA(w)}{\pi}\]
where $\|\sum a_jz^j\|_2^2=\sum|a_j|^2$ and $dA(w)$ is Lebesgue area measure on the unit disk.
\end{theorem}

The following \textbf{Littlewood Subordination Theorem} can be used to show that composition operators are bounded in a variety of spaces.
\begin{theorem}[Theorem 3.8, \cite{Cowen}]\label{Littlewood}
Let $\varphi$ be an inner function on the disk, and $C_\varphi$ be the composition operator defined by $C_\varphi f=f\circ \varphi, f\in Hol(\mathbb{D})$. If $f\in H^p(\mathbb{D})$, $0<p<\infty$, then
\begin{equation*}
       \left(\frac{1-|\varphi(0)|}{1+|\varphi(0)|}\right)^{1/p}\|f\|_p   \leq\|C_\varphi f\|_p\leq \left(\frac{1+|\varphi(0)|}{1-|\varphi(0)|}\right)^{1/p}\|f\|_p
\end{equation*}
In particular, if  $\varphi(0)=0$, then $C_\varphi$ is an isometry.
\end{theorem}

\subsection{Core operator and core function}
\begin{definition}\label{def2.2}
Let $M$ be a submodule in $H^2(\mathbb{D}^2)$ and $K^M_{\lambda,\mu}(z,w)$ be the reproducing kernel for $M$. Then the core function $G^M_{\lambda,\mu}(z,w)$ for $M$ is defined as
\begin{equation}\label{2.1}
  G^M_{\lambda,\mu}(z,w):=(1-\bar{\lambda}z)(1-\bar{\mu}w)K^M_{\lambda,\mu}(z,w).
\end{equation}
With the core function $G^M_{\lambda,\mu}(z,w)$, the associated core operator $C_M$ is defined by
\begin{equation}\label{2.2}
  C_M(f)(z):=\int_{\mathbb{T}^2} G^M_{\lambda,\mu}(z,w) f(\lambda,\mu) |d\lambda||d\mu|,~~~f\in H^2(\mathbb{D}^2).
\end{equation}

\end{definition}

The core function and core operator are introduced by K. Guo and R. Yang \cite{Core} and have been well studied in  \cite{Core,Ya2004,Ya2005}.  It is known that for every submodule $M$ in $H^2(\mathbb{D}^2)$, $G^M_{z,w}(z,w)=1$ a.e. on $\mathbb{T}^2$, and the core operator $C_M$ is zero on $M^\perp=H^2(\mathbb{D}^2)\ominus M$. If we denote by $R_z, R_w$ and $S_z,S_w$ the restriction of the multiplication operators $M_z,M_w$ on $M$ and respectively the compression of $M_z,M_w$ to $M^\perp$, then on $M$,
\[C_M=I-R_zR_z^*-R_wR_w^*+R_zR_wR_z^*R_w^*.\]
 In particular, it was shown \cite{Ya2005} that on $M$, $C_M^2$ is unitarily equivalent to
\begin{equation*}
  \left(
    \begin{array}{cc}
      [R_z^*, R_z][R_w^*,R_w][R_z^*,R_z] & 0 \\
      0 & [R_z^*,R_w][R_w^*,R_z] \\
    \end{array}
  \right),
\end{equation*}
here $[A,B]:=AB-BA$ denotes the commutator of $A,B$. Therefore $C_M$ is Hilbert-Schmidt if and only if $[R_z^*, R_z][R_w^*,R_w]$ and $[R_z^*,R_w]$ are both Hilbert-Schmidt. Thus $M$ is a Hilbert-Schmidt if and only if the core operator $C_M$ is Hilbert-Schmidt. Moreover, in this case,
\begin{equation}\label{2.3}
  \|C_M\|_{H.S.}^2=\int_{\mathbb{T}^2}\int_{\mathbb{T}^2} ~|G^M_{\lambda,\mu}(z,w)|^2 ~|d\lambda||d\mu||dz||dw|=\Sigma_0^M+\Sigma_1^M . 
\end{equation}
If $C_M$ is trace class, then 
\begin{equation}\label{2.4}
 \tr(C_M)=\Sigma_0^M-\Sigma_1^M=\lim_{r\rightarrow 1}\int_{\mathbb{T}^2} G^M_{rz,rw}(rz,rw) |dz||dw|=1.
\end{equation}

\subsection{Reproducing kernel Hilbert spaces}
In this subsection, we will present some basic theory of reproducing kernel Hilbert spaces.
For more information, see \cite{Aron, Pa}.

Let $\Omega \subset \mathbb{C}^d$. We say a function $K: \Omega\times \Omega\rightarrow \mathbb{C}$ is a positive kernel on $\Omega$ if it is self-adjoint ($K(x,y)=\overline{K(y,x)}$), and for any finite sets $\{\lambda_1, \lambda_2,\cdots,\lambda_m\}\subset \Omega$, the matrix $(K(\lambda_i,\lambda_j))_{i,j=1}^m$ is positive semi-definite. For convenience, we write this property as $K\geq 0$. If $K_1,K_2$ are two positive kernel on $\Omega$, $K_1\geq K_2$ means that $K_1-K_2$ is a positive kernel.

There are some general methods to construct new positive kernels from old ones.
\begin{lemma}
Let $K_1,K_2$ be two positive kernel on $\Omega$. Then
\begin{enumerate}
  \item $K_1+K_2\geq0$, $K_1 \cdot K_2 \geq 0$.
  \item If $f:\Omega \rightarrow \mathbb{C}$ is a function, then $\overline{f(y)}K(x,y)f(x) \geq 0$.
  \item If $B:\Omega \rightarrow \Omega $ is a self-map, then $K\circ B \geq 0$, here $K\circ B(x,y)=K(B(x),B(y))$.
\end{enumerate}
\end{lemma}
A reproducing kernel Hilbert space $\mathcal{H}$ on $\Omega$ is a Hilbert space of complex valued functions on $\Omega$ such that every point evaluation is a continuous linear functional. Thus by Riesz representation theorem, for each point $y\in \Omega$, there exists an $K_y\in \mathcal{H}$ such that for each $f\in \mathcal{H}$,
\[ \langle  f,K_y \rangle_{\mathcal{H}}=f(y). \]
Since $K_y(x)=\langle K_y, K_x \rangle_{\mathcal{H}}$, $K$ can be regarded as a function on $\Omega\times \Omega$ and we write $K_y(x)$ as $K(x,y)$. Such $K$ is a positive kernel and the Hilbert space $\mathcal{H}$ with reproducing kernel $K$ is denoted by $\mathcal{H}(K)$. The following theorem, due to Moore \cite{Moore}, shows that there is one-to-one correspondence between reproducing kernel Hilbert spaces and positive kernels.
\begin{theorem}\label{unique}
Let $K: \Omega\times \Omega\rightarrow \mathbb{C}$ be a positive kernel. Then there exists a unique reproducing kernel Hilbert space $\mathcal{H}(K)$ whose reproducing kernel is $K$.
\end{theorem}


\begin{definition}
Let $\mathcal{H}$ be a  Hilbert space  with inner product $\langle \cdot,\cdot\rangle$. A set of vectors $\{ f_s: ~s\in S \}\subseteq \mathcal{H}$ is called a \textbf{Parseval frame} for $\mathcal{H}$ provided that
\[ \|h\|^2= \sum_{s\in S} |\langle h,f_s\rangle|^2  \]
for every $h\in \mathcal{H}$.
\end{definition}
It is easy to see that $\|f_s\|\leq 1$ and if every $f_s$ has unit norm, then $\{f_s: ~s\in S \}$ is an orthonormal basis of $\mathcal{H}$.

\begin{theorem}[Papadakis, Theorem 2.10, \cite{Pa}]\label{PF}
Let $\mathcal{H}(K)$ be a reproducing kernel Hilbert space and let $\{ f_s: ~s\in S \}\subseteq \mathcal{H}$. Then $\{ f_s: ~s\in S \}$ is a Parseval frame for $\mathcal{H}$ if and only if
\[ K(x,y)=\sum_{s\in S} \overline{f_s(y)} f_s(x) ,  \]
where the series converges pointwise.
\end{theorem}
In fact, the hypothesis that each $f_s\in \mathcal{H}$ is redundant. The following theorem shows that $K(x,y)=\sum_{s\in S}  \overline{f_s(y)}f_s(x)$ automatically implies $f_s\in \mathcal{H}$.

\begin{theorem}[Theorem 3.11, \cite{Pa}]\label{belong}
Let $\mathcal{H}(K)$ be a reproducing kernel Hilbert space on
$\Omega$ and let $f:\Omega \rightarrow \mathbb{C}$ be a function. Then $f\in \mathcal{H}(K)$ with $\|f\|_{\mathcal{H}}\leq c$ if and only if
\[ c^2K(x,y)- \overline{f(y)} f(x) \geq 0. \]
\end{theorem}

\begin{definition}
Let $\mathcal{H}_i, i=1,2$ be reproducing kernel Hilbert spaces on the same set $\Omega$ and let $K_i,i=1,2$ denote their kernel functions. A function $f:\Omega\rightarrow \mathbb{C}$ is called a \textbf{multiplier} of $\mathcal{H}_1$ into $\mathcal{H}_2$ provided that $f\mathcal{H}_1:=\{fh:h\in \mathcal{H}_1\}\subseteq \mathcal{H}_2$.
\end{definition}

 We let $\mathcal{M}(\mathcal{H}_1,\mathcal{H}_2)$ denote the set of all multipliers of $\mathcal{H}_1$ into $\mathcal{H}_2$. When $\mathcal{H}_1=\mathcal{H}_2=\mathcal{H}$ with reproducing kernel $K_1=K_2=K$, then we call a multiplier of $\mathcal{H}$ into $\mathcal{H}$, more simply, a multiplier of $\mathcal{H}$ and denote $\mathcal{M}(\mathcal{H},\mathcal{H})$ by $\mathcal{M}(\mathcal{H})$. Given a multiplier $f\in \mathcal{M}(\mathcal{H}_1,\mathcal{H}_2)$, let $M_f:\mathcal{H}_1\rightarrow \mathcal{H}_2$ denote the multiplication operator $M_f(h)=fh$.  Clearly, the set of multipliers $\mathcal{M}(\mathcal{H}_1,\mathcal{H}_2)$ is a vector space and the set of multipliers $\mathcal{M}(\mathcal{H})$ is an algebra. The following theorem gives a characterization of multipliers.
\begin{theorem}[Theorem 5.21, \cite{Pa}]\label{Mult}
Let $\mathcal{H}_i, i=1,2$ be reproducing kernel Hilbert spaces on  $\Omega$ with kernels $K_i,i=1,2$ and let $f:\Omega\rightarrow \mathbb{C}$ be a function. The following are equivalent:
\begin{enumerate}
  \item $f\in \mathcal{M}(\mathcal{H}_1,\mathcal{H}_2)$;
  \item $f\in \mathcal{M}(\mathcal{H}_1,\mathcal{H}_2)$ and $M_f$ is a bounded operator;
  \item there is a constant $c\geq 0$, such that $\overline{f(y)} K_1(x,y) f(x)\leq c^2K_2(x,y)$.
\end{enumerate}
Moreover, in this case, $\|M_f\|$ is the least constant $c$ satisfying the inequality in $(3)$.
\end{theorem}

\begin{corollary}\label{Contain}
Let $\mathcal{H}(K_1)$ and $\mathcal{H}(K_2)$ be two reproducing kernel Hilbert spaces on $\Omega$. Then $\mathcal{H}(K_1)\subset \mathcal{H}(K_2)$ if and only if there is some constant $c\geq 0$ such that
\[ K_1(x,y)\leq c^2 K_2(x,y). \]
Moreover, in this case, for any $h\in \mathcal{H}(K_1)$, $\|h\|_{\mathcal{H}(K_2)}\leq c \|h\|_{\mathcal{H}(K_1)}$.
\end{corollary}

\section{The  reproducing kernel of submodule $M_{\theta,\varphi}$}
Let $\theta(z),\varphi(w)$ be two nonconstant inner functions and $M$ be a submodule in $H^2(\mathbb{D}^2)$. Let $M_{\theta,\varphi}$ denote the submodule $[C_{\theta,\varphi}M]$ (by Theorem\ref{Littlewood}, $C_{\theta,\varphi}$ is bounded, hence $C_{\theta,\varphi}M \subset H^2(\mathbb{D}^2)$), that is, the smallest submodule that contains $C_{\theta,\varphi}M$.
\begin{theorem}\label{thm3.1}
 Let $K^M_{\lambda,\mu}(z,w)$ and $K^{M_{\theta,\varphi}}_{\lambda,\mu}(z,w)$ be the reproducing kernel for $M$ and  $M_{\theta,\varphi}$, respectively. Then
\begin{equation}\label{3.1}
  K^{M_{\theta,\varphi}}_{\lambda,\mu}(z,w)  =\frac{1-\overline{\theta(\lambda)}\theta(z)}{1-\bar{\lambda}z} \frac{1-\overline{\varphi(\mu)}\varphi(w)}{1-\bar{\mu}w} K^{M}_{\theta(\lambda),\varphi(\mu)}(\theta(z),\varphi(w)).
\end{equation}

\end{theorem}
\begin{remark}
For simplicity, we set $R_{\lambda,\mu}(z,w)=\frac{1-\overline{\theta(\lambda)}\theta(z)}{1-\bar{\lambda}z} \frac{1-\overline{\varphi(\mu)}\varphi(w)}{1-\bar{\mu}w}$ (clearly, $R\geq 0$) and let $B=(\theta,\varphi)$ be the symbol of composition operator $C_{\theta,\varphi}$, then the main result of above theorem can be written simply as
\[ K^{M_{\theta,\varphi}}= K^M \circ B~ \cdot R. \]
\end{remark}
\begin{proof}
Let $S_{\lambda,\mu}(z,w):=\frac{1}{(1-\bar{\lambda}z)(1-\overline{\mu}w)}$ be the Szeg\"{o} kernel of $H^2(\mathbb{D}^2)$, as $M$ is a submodule of $H^2(\mathbb{D}^2)$, we have $K^M \circ B\leq S\circ B$. Note that $R=\frac{S}{S\circ B}$, hence
\[ K^M \circ B~ \cdot R \leq S\circ B ~\cdot R =S.\]
By Corollary \ref{Contain}, it implies that $\mathcal{H}(K^M \circ B~ \cdot R)$  is contractively contained in $H^2(\mathbb{D}^2)$, that is,  $\mathcal{H}(K^M \circ B~ \cdot R)\subseteq H^2(\mathbb{D}^2)$ and for $h\in \mathcal{H}(K^M \circ B~ \cdot R)$,
\begin{equation}\label{3.2}
  \|h\|_2\leq \|h\|_{\mathcal{H}(K^M \circ B~ \cdot R)}.
\end{equation}

Let $\{\alpha_m(z)\}$, $\{\beta_n(w)\}$ and $\{e_k(z,w)\}$ be orthonormal bases for $\mathcal{K}_\theta:=H^2_z \ominus \theta H^2_z$, $\mathcal{K}_\varphi:=H^2_w \ominus \varphi H^2_w$ and $M$, respectively. By  Theorem\ref{PF},
$\{ \alpha_m(z) \beta_n(w) e_k(\theta(z),\varphi(w)) : m,n,k \geq 0 \}$ is a Parseval frame of $\mathcal{H}(K^M \circ B~ \cdot R)$. Write $e_k(z,w)=\sum_{i,j\geq 0} a^k_{ij} z^iw^j$, by a direct computation, we can verify that
\[ \langle \alpha_m \beta_n e_k(\theta,\varphi),\alpha_{m'} \beta_{n'} e_{k'}(\theta,\varphi)\rangle_{H^2(\mathbb{D}^2)}=\delta_{m,n,k}^{m',n',k'}, \]
where $\delta_{m,n,k}^{m',n',k'} \neq 0$ if and only if $m=m',n=n',k=k'$ and $\delta_{m,n,k}^{m,n,k}=1$.
That is, $\{ \alpha_m(z) \beta_n(w) e_k(\theta(z),\varphi(w)) : m,n,k \geq 0 \}$ is an orthonomal set in $H^2(\mathbb{D}^2)$.

Combine the fact that $\{ \alpha_m(z) \beta_n(w) e_k(\theta(z),\varphi(w)) : m,n,k \geq 0 \}$ is a Parseval frame of $\mathcal{H}(K^M \circ B~ \cdot R)$ and formula $(\ref{3.2})$, we get that $\{ \alpha_m(z) \beta_n(w) e_k(\theta(z),\varphi(w)) : m,n,k \geq 0 \}$ is also an orthonormal basis of $\mathcal{H}(K^M \circ B~ \cdot R)$. This fact implies that $\mathcal{H}(K^M \circ B~ \cdot R)$ is exactly a closed subspace of $H^2(\mathbb{D}^2)$, not only just as a set, but also inheriting the Hilbert space structure of $H^2(\mathbb{D}^2)$. In other words, $\mathcal{H}(K^M \circ B~ \cdot R)$ is isometrically contained in $H^2(\mathbb{D}^2)$. Therefore by Theorem \ref{unique}, our conclusion is equivalent to that $\mathcal{H}(K^M \circ B~ \cdot R)$ equals $M_{\theta,\varphi}$ as a set.

Since $\frac{1-\overline{\varphi(\mu)}\varphi(w)}{1-\bar{\mu}w}$ and $K^{M\ominus zM}$ are both positive kernels, we get
\begin{align*}
  (1-\bar{\lambda}z) K^M \circ B~ \cdot R &=\frac{1-\overline{\varphi(\mu)}\varphi(w)}{1-\bar{\mu}w} \cdot (1-\overline{\theta(\lambda)}\theta(z))K^M \circ B\\
   &=\frac{1-\overline{\varphi(\mu)}\varphi(w)}{1-\bar{\mu}w} \cdot K^{M\ominus zM} \circ B \geq 0.
\end{align*}
Then by Theorem \ref{Mult}, it implies that $\mathcal{H}(K^M \circ B~ \cdot R)$ is $z$-invariant. Similarly, we can also show that $\mathcal{H}(K^M \circ B~ \cdot R)$ is $w$-invariant.  For any fixed point $(\lambda,\mu)\in \mathbb{D}^2$, note that $R_{\lambda,\mu}$ is invertible in $H^\infty(\mathbb{D}^2)$, so there is a polynomial sequence $\{p_n\}$ with $\|p_n\|_\infty\leq \|R^{-1}_{\lambda,\mu}\|_\infty$ such that $p_n$ converges to $R^{-1}_{\lambda,\mu}$ in $2$-norm. For any fixed $f\in H^2(\mathbb{D}^2)$ and any $\varepsilon >0$, there is polynomial $q$ such that $\|q-f\|_2< \varepsilon$, and then
\begin{align*}
  \|p_nf-R^{-1}_{\lambda,\mu} f \|_2 &\leq \|p_nf-p_nq\|_2+\|p_nq-R^{-1}_{\lambda,\mu}q\|_2+\|R^{-1}_{\lambda,\mu}q-R^{-1}_{\lambda,\mu}f\|_2 \\
   &\leq \|R^{-1}_{\lambda,\mu}\|_\infty \varepsilon +\|q\|_\infty \|p_n-R^{-1}_{\lambda,\mu}\|_2+\|R^{-1}_{\lambda,\mu}\|_\infty \varepsilon\longrightarrow 2\|R^{-1}_{\lambda,\mu}\|_\infty \varepsilon.
\end{align*}
This implies that $T_{p_n}$ converges to $T_{R^{-1}_{\lambda,\mu}}$ in strong operator topology. Thus we have
\[  p_n R_{\lambda,\mu}\cdot (K^M\circ B)_{\lambda,\mu} \xrightarrow{\|\cdot\|_2}  K^M_{B(\lambda,\mu)}(B(z,w))=C_B K^M_{B(\lambda,\mu)}, \]
where $C_B$ is the composition operator. Since $\mathcal{H}(K^M \circ B~ \cdot R)$ is $z,w$-invariant, we obtain that $C_B K^M_{B(\lambda,\mu)}\in \mathcal{H}(K^M \circ B~ \cdot R)$ for any $(\lambda,\mu)\in \mathbb{D}^2$. As $B$ has dense range in $\mathbb{D}^2$, we get that $C_BM\subseteq \mathcal{H}(K^M \circ B~ \cdot R)$. Again by the fact that $\mathcal{H}(K^M \circ B~ \cdot R)$ is $z,w$-invariant, we deduce that $M_{\theta,\varphi}\subseteq \mathcal{H}(K^M \circ B~ \cdot R)$.

On the other hand, for any point $(\lambda,\mu)\in \mathbb{D}^2$, similar arguments show that there is a polynomial sequence $\{p_n\}$ with $\|p_n\|_\infty\leq \|R_{\lambda,\mu}\|_\infty$ such that $p_n$ converges to $R_{\lambda,\mu}$ in $2$-norm. Thus we have
\[  p_n \cdot (K^M\circ B)_{\lambda,\mu} \xrightarrow{\|\cdot\|_2} R_{\lambda,\mu} (K^M\circ B)_{\lambda,\mu}.  \]
As $(K^M\circ B)_{\lambda,\mu}=C_B K^M_{B(\lambda,\mu)}\in M_{\theta,\varphi}$, we get that $R_{\lambda,\mu} \cdot (K^M\circ B)_{\lambda,\mu}$ is contained in the submodule $M_{\theta,\varphi}$ for any $(\lambda,\mu)\in \mathbb{D}^2$. It implies that $\mathcal{H}(K^M \circ B~ \cdot R)\subseteq M_{\theta,\varphi}$.
\end{proof}

\begin{remark}
Let $f$ be in $\mathcal{H}(R)=\mathcal{K_\theta}\otimes \mathcal{K}_\varphi$ with unit norm, then the weighted composition operator $M_fC_B$ is an  isometry from $M$ into $M_{\theta,\varphi}$. Indeed, by Theorem \ref{belong}, $R-\overline{f(\lambda,\mu)}f(z,w)\geq 0$, so
\[ \overline{f(\lambda,\mu)}f(z,w) K^M\circ B \leq K^M\circ B \cdot R,  \]
which implies that $f$ is a contractive multiplier of $\mathcal{H}(K^M\circ B)$ into $M_{\theta,\varphi}$ (Theorem \ref{Mult}). Therefore, it suffices to show that $M_fC_B$ is an isometry on $H^2(\mathbb{D}^2)$. Take $M=H^2(\mathbb{D}^2)$, the above inequality show that $M_fC_B$ is a contraction, so it is sufficient to show that $\|M_fC_Bp\|=\|p\|$ for every polynomial $p$. Writing $f=\sum_{m,n\geq 0}f_{mn}\alpha_m(z)\beta_n(w)$, by direct computation, we can check that $\|M_fC_Bp\|^2=\|p\|^2$ holds.\\
\end{remark}

Due to Theorem \ref{thm3.1}, we naturally have the following fact.
\begin{corollary}\label{Cor3.1}
Let $G^M$ and  $G^{M_{\theta,\varphi}}$ be the core function of submodule $M$ and $M_{\theta,\varphi}$, respectively. Then
\begin{equation}\label{3.3}
  G^{M_{\theta,\varphi}}_{\lambda,\mu}(z,w)=G^{M}_{\theta(\lambda),\varphi(\mu)}(\theta(z),\varphi(w)),
\end{equation}
or simply, $G^{M_{\theta,\varphi}}=G^M\circ B$.
\end{corollary}

Moreover, by the Littlewood's Subordination Principle, we easily deduce the following result.
\begin{theorem}\label{thm3.3}
Let $\theta(z), \varphi(w)$ be two nonconstant inner functions and $M$ be a submodule in $H^2(\mathbb{D}^2)$. Then $M$ is Hilbert-Schmidt if and only if $M_{\theta,\varphi}$ is Hilbert-Schmidt. Moreover,
\begin{equation}\label{3.4}
 \frac{1-|\theta(0)|}{1+|\theta(0)|}\cdot \frac{1-|\varphi(0)|}{1+|\varphi(0)|} \cdot \|C_{M}\|_{H.S.} \leq \|C_{M_{\theta,\varphi}}\|_{H.S.}\leq \frac{1+|\theta(0)|}{1-|\theta(0)|}\cdot \frac{1+|\varphi(0)|}{1-|\varphi(0)|} \cdot \|C_{M}\|_{H.S.}.
\end{equation}

\end{theorem}
\begin{proof}
It is immediately obtained from the Theorem \ref{Littlewood} and the equations \eqref{2.3}, \eqref{3.3}.
\end{proof}

\begin{corollary}
Let $\theta(z), \varphi(w)$ be two nonconstant inner functions and $q$ be a polynomial, then $[q(\theta,\varphi)]$ is a Hilbert-Schmidt submodule.
\end{corollary}
\begin{proof}
It is immediately obtained from Theorem \ref{thm3.3} and the fact that every polynomial submodule is Hilbert-Schmidt (see Corollary 2.5 in \cite{Ya3}).
\end{proof}

\section{Numerical invariant functions}
In this section, we always assume the submodule $M$ is Hilbert-Schmidt. In this case,
\[ [R_z^*,R_z][R_w^*,R_w], ~~[R_z^*,R_w] \]
are both Hilbert-Schmidt. Let $\varphi_a(z):=\frac{a-z}{1-\overline{a}z}$, then it is not hard to check that the perturbations $$[R_{\varphi_a(z)}^*,R_{\varphi_a(z)}][R_{\varphi_b(w)}^*,R_{\varphi_b(w)}],~~[R_{\varphi_a(z)}^*, R_{\varphi_b(w)}]$$
are also Hilbert-Schmidt for any $(a,b)\in \mathbb{D}^2$.

For $(a,b) \in \mathbb{D}^2$, let $P_a$ and $Q_b$ denote the orthogonal projections from $M$ onto $M\ominus (z-a)M$ and $M\ominus (w-b)M$, respectively. Clearly,
\[
P_a=[R_{\varphi_a(z)}^*, R_{\varphi_a(z)}],~~ ~Q_b=[R_{\varphi_b(w)}^*, R_{\varphi_b(w)}].
\]
For $(a,b) \in \mathbb{D}^2$, let $P'_a$ and $Q'_b$ denote the orthogonal projections from $M$ onto $(z-a)(M\ominus (w-b)M)$ and $(w-b)(M\ominus (z-a)M)$, respectively. Similarly, we have
\[
P'_a=R_{\varphi_a(z)}Q_bR_{\varphi_a(z)}^*,~~ ~Q'_b=R_{\varphi_b(w)}P_aR_{\varphi_b(w)}^*.
\]

\begin{definition}
For a Hilbert-Schmidt submodule $M$, we define the \textbf{numerical invariant functions} $\Sigma_0^M(a,b)$, $\Sigma_1^M(a,b)$ as
\[  \Sigma^M_0(a,b )=\|P_a Q_b\|_{H.S.}^2=\|[R_{\varphi_a(z)}^*, R_{\varphi_a(z)}][R_{\varphi_b(w)}^*,R_{\varphi_b(w)}] \|_{H.S.}^2 ,\]
\[\Sigma^M_1(a,b )=\|P'_a Q'_b\|_{H.S.}^2=\|[R_{\varphi_a(z)}^*, R_{\varphi_b(w)}]\|_{H.S.}^2,~~~(a,b)\in \mathbb{D}^2.\]
When $(a,b)=(0,0)$, they are exactly the numerical invariants first introduced by R. Yang \cite{Ya3}. Moreover, we can also define the \textbf{core operator function} as
\[C_M(a,b)=I-R_{\varphi_a(z)}R_{\varphi_a(z)}^*-R_{\varphi_b(w)}R_{\varphi_b(w)}^*
+R_{\varphi_a(z)}R_{\varphi_b(w)}R_{\varphi_a(z)}^*R_{\varphi_b(w)}^* .\]
When $(a,b)=(0,0)$, it is exactly the core operator $C_M$ for submodule $M$. 
\end{definition}
By some computation, it is not hard to check that $ C_M(a,b)=C_{\varphi_a,\varphi_b}C_MC^*_{\varphi_a,\varphi_b}$
and $C_M^2(a,b)$ is unitarily equivalent to 
\begin{equation*}
  \left(
    \begin{array}{cc}
      [R_{\varphi_a(z)}^*, R_{\varphi_a(z)}][R_{\varphi_b(w)}^*,R_{\varphi_b(w)}][R_{\varphi_a(z)}^*,R_{\varphi_a(z)}] & 0 \\
      0 & [R_{\varphi_a(z)}^*,R_{\varphi_b(w)}][R_{\varphi_b(w)}^*,R_{\varphi_a(z)}] \\
    \end{array}
  \right).
\end{equation*}
Thus we still have 
\[ \|C_M(a,b)\|_{H.S.}^2=\Sigma^M_0(a,b )+\Sigma^M_1(a,b ). \]
Further, if $C_M$ is trace class, then $C_M(a,b)$ is also trace class for any $(a,b)\in \mathbb{D}^2$ and
\begin{equation*}
 \tr(C_M(a,b))=\Sigma_0^M(a,b)-\Sigma_1^M(a,b).
\end{equation*}


In \cite{Luo}, the authors generalized the definition of fringe operator which was first introduced in \cite{Ya3}. For $(a,b)\in \mathbb{D}^2$, the fringe operator $F_{a,b}$ on $M\ominus (z-a)M$ is defined by
\[ F_{a,b}f=P_a R_{\varphi_b(w)} f,~~~~~f\in M\ominus (z-a)M.  \]

It was shown \cite{Luo} that $F_{a,b}$  is Fredholm if and only if the pair $(R_{\varphi_a(z)},R_{\varphi_b(w)})$ is Fredholm, and in this case
\begin{equation}\label{trace}
  \tr [F_{a,b}^*, F_{a,b}]=-\ind F_{a,b}=\ind (R_{\varphi_a(z)},R_{\varphi_b(w)})= \Sigma^M_0(a,b )- \Sigma^M_1(a,b ).
\end{equation}

For the computation of $\Sigma^M_0(a,b ),\Sigma^M_1(a,b )$, we need a simple fact.
\begin{lemma}\label{lemma4.4}
Let $\mathcal{H}$ be a reproducing Hilbert subspace with reproducing kernel $K$ and let $P,Q$ be two orthogonal projections onto closed subspaces $M,N\subset \mathcal{H}$, respectively.  Suppose $\{e_k:k\geq 0\}$, $\{f_l:l\geq 0\}$ are  Parseval frames for $M,N$, respectively. If $PQ$ is Hilbert-Schmidt, then
\[\|PQ\|_{H.S.}^2 =\sum_{k,l\geq 0} |\langle e_k, f_l \rangle|^2.  \]
\end{lemma}
\begin{proof}
Let $\{g_m:m\geq 0\}$ be an orthonormal basis of $N$, then we have
\begin{align*}
  \|PQ\|_{H.S.}^2  = \sum_{m\geq 0} \|Pg_m\|^2=\sum_{m,k\geq 0} |\langle g_m, e_k \rangle|^2=\sum_{k\geq 0}\|Qe_k\|^2=\sum_{k,l\geq 0} |\langle e_k, f_l \rangle|^2.
\end{align*}
\end{proof}

\begin{proposition}\label{prop4.5}
For a Hilbert-Schmidt submodule $M$, suppose $\{e_k:k\geq 0\}$, $\{f_l:l\geq 0\}$ are Parseval frames for $M\ominus zM, M\ominus wM$, respectively. Then
\begin{equation}\label{4.5}
  \Sigma_0^M(a,b)=(1-|a|^2)(1-|b|^2) \sum_{k,l\geq 0} |\langle \frac{e_k}{1-\bar{a}z}, \frac{f_l}{1-\bar{b}w} \rangle|^2,
\end{equation}
\begin{equation}\label{4.6}
  \Sigma_1^M(a,b)=(1-|a|^2)(1-|b|^2) \sum_{k,l\geq 0} |\langle \frac{w e_k}{1-\bar{b}w}, \frac{zf_l}{1-\bar{a}z} \rangle|^2.
\end{equation}
\end{proposition}
\begin{proof}
Set $\alpha(z)=\frac{\sqrt{1-|a|^2}}{1-\bar{a}z}, \beta(w)=\frac{\sqrt{1-|b|^2}}{1-\bar{b}w}$, then
\[ K^{M\ominus (z-a)M} =(1-\overline{\varphi_a(\lambda)}\varphi_a(z))K^M=
\frac{(1-|a|^2)(1-\bar{\lambda}z)}{(1-a\bar{\lambda})(1-\bar{a}z)}K^{M}
=\overline{\alpha(\lambda)}\alpha(z)K^{M\ominus zM},  \]
\[ K^{M\ominus (w-b)M} =(1-\overline{\varphi_b(\mu)}\varphi_b(w))K^M=
\frac{(1-|b|^2)(1-\bar{\mu}w)}{(1-b\bar{\mu})(1-\bar{b}w)}K^{M}
=\overline{\beta(\mu)}\beta(w)K^{M\ominus wM}.  \]
So by Theorem \ref{PF} and Lemma \ref{lemma4.4}, we get the formula \eqref{4.5}. Similarly,
\[ K^{(w-b)(M\ominus(z-a)M)}=
\overline{\varphi_b(\mu)}\varphi_b(w)\overline{\alpha(\lambda)}\alpha(z)K^{M\ominus zM}, \]
\[ K^{(z-a)(M\ominus(w-b)M)}=
\overline{\varphi_a(\lambda)}\varphi_a(z)\overline{\beta(\mu)}\beta(w)K^{M\ominus wM}. \]
Thus
\begin{align*}
   \Sigma_1^M(a,b)&=(1-|a|^2)(1-|b|^2) \sum_{k,l\geq 0} |\langle \frac{\varphi_b(w) e_k}{1-\bar{a}z}, \frac{\varphi_a(z)f_l}{1-\bar{b}w} \rangle|^2\\
   & =(1-|a|^2)(1-|b|^2) \sum_{k,l\geq 0} |\langle \frac{ e_k}{a-z}, \frac{f_l}{b-w} \rangle_{L^2(\mathbb{T}^2)}|^2\\
   & =(1-|a|^2)(1-|b|^2) \sum_{k,l\geq 0} |\langle \frac{ \bar{z}e_k}{a\bar{z}-1}, \frac{\bar{w}f_l}{b\bar{w}-1} \rangle_{L^2(\mathbb{T}^2)}|^2\\
   &=(1-|a|^2)(1-|b|^2) \sum_{k,l\geq 0} |\langle \frac{w e_k}{1-\bar{b}w}, \frac{zf_l}{1-\bar{a}z} \rangle|^2.
\end{align*}
\end{proof}

\begin{remark}
In \cite{Ya3}, R. Yang also define the higher-order numerical invariants 
\[ \Sigma_k=\sum_{i,j=0}^{\infty} |\langle w^k e_i, z^k f_j\rangle|^2,~~~~~~k\geq 2, \]
where $\{e_i : i\geq 0\}$, $\{f_j : j\geq 0\}$ are Parseval frames for $M\ominus zM, M\ominus wM$, respectively. Based on concrete computations and observations of some typical submodules, R. Yang proposed an attractive but seemingly "simple" conjecture:
\[ \{\Sigma_k(M): k\geq 0\}~\text{is a decreasing sequence for every submodule}~M. \]
However, due to the difficulty of computing higher-order numerical invariants for general submodules and the lack of effective tools and methods, these higher-order invariants have remained largely unexplored over the past two decades. As a result, this conjecture remains an open problem.

In fact, the numerical invariant functions can also be naturally extended to the higher-order cases. Let $k_a(z):=\frac{\sqrt{1-|a|^2}}{1-\bar{a}z}$ be the normalized Szegö kernel of $H^2(\mathbb{D})$, then the formulas \eqref{4.5} and \eqref{4.6} can be write as 
\begin{align*}
\Sigma_0^M(a,b)&= \sum_{i,j\geq 0} |\langle k_a(z) e_i, k_b(w) f_j \rangle|^2\\
\Sigma_1^M(a,b)&=\sum_{i,j\geq 0} |\langle wk_b(w) e_i, zk_a(z) f_j \rangle|^2.
\end{align*}
Therefore, we can similarly define higher-order numerical invariant functions
\[ \Sigma_k^M(a,b)=\sum_{i,j\geq 0} |\langle w^k k_b(w) e_i, z^k k_a(z) f_j \rangle|^2,~~~~~~k\geq 2.  \]
Compared to numerical invariants, numerical invariant functions allow us to study the properties of submodules more thoroughly from Fourier analysis and functional theory perspectives. 
\end{remark}

\begin{example}

\begin{enumerate}
\item Consider the submodule $M=zH^2(\mathbb{D}^2)+wH^2(\mathbb{D}^2)$, by a direct computation, one can verify that
\[\Sigma_0^M(a,b)=(1-|a|^2)(1-|b|^2)+1,~~~\Sigma_1^M(a,b)=(1-|a|^2)(1-|b|^2).\]
\item Consider the Beurling type submodule $M=\theta H^2(\mathbb{D}^2)$, it is easy to check that $\Sigma_0^M(a,b)\equiv 1$ and $\Sigma_0^M(a,b)\equiv 0$.
\end{enumerate}
\end{example}

It is known that the classical numerical invariants $\Sigma_0^M,\Sigma_1^M$ satisfy the equations (see Theorem 4.4 in \cite{Ya2004})
\[ \Sigma_0^M-\Sigma_1^M=1=\ind(R_z,R_w).\]
So it a natural question whether similar equations still hold for the numerical invariant functions $\Sigma_0^{M}(a,b)$ and $\Sigma_1^{M}(a,b)$.

\begin{theorem}\label{thm4.6}
Let  $\theta,\varphi$ be two nonconstant one-variable inner functions. If $M$ is Hilbert-Schmidt submodule, then
\[\Sigma_0^{M_{\theta,\varphi}}(0,0)=\Sigma_0^M(\theta(0),\varphi(0)),~~~
\Sigma_1^{M_{\theta,\varphi}}(0,0)=\Sigma_1^M(\theta(0),\varphi(0)). \]
Moreover, $\|C_{M_{\theta,\varphi}}\|_{H.S.}^2=\Sigma_0^M(\theta(0),\varphi(0))+\Sigma_1^M(\theta(0),\varphi(0)).$
\end{theorem}
\begin{proof}

As usual, set $B(z,w)=(\theta(z),\varphi(w))$. By Theorem \ref{thm3.1},
$$K^{M_{\theta,\varphi}}_{\lambda,\mu}(z,w)  =\frac{1-\overline{\theta(\lambda)}\theta(z)}{1-\bar{\lambda}z} \frac{1-\overline{\varphi(\mu)}\varphi(w)}{1-\bar{\mu}w} K^{M}_{\theta(\lambda),\varphi(\mu)}(\theta(z),\varphi(w)),$$
 so we have
\[ K^{M_{\theta,\varphi}\ominus zM_{\theta,\varphi}}= K^{M\ominus zM}\circ B \cdot \frac{1-\overline{\varphi(\mu)}\varphi(w)}{1-\bar{\mu}w}, \]
\[ K^{M_{\theta,\varphi}\ominus wM_{\theta,\varphi}}= K^{M\ominus wM}\circ B \cdot \frac{1-\overline{\theta(\lambda)}\theta(z)}{1-\bar{\lambda}z}. \]

Let $\{\alpha_m(z)\}$, $\{\beta_n(w)\}$ , $\{e_k(z,w)\}$  and $\{f_l(z,w)\}$ be orthonormal bases of $\mathcal{K}_\theta$, $\mathcal{K}_\varphi$, $M\ominus zM$ and $M\ominus wM$, respectively. Then by Theorem \ref{PF}, $\{ \beta_n(w) e_k(\theta,\varphi): n,k\geq 0\}$ and $\{\alpha_m(z) f_l(\theta,\varphi): m,l \geq 0\}$ are Parseval frames for $M_{\theta,\varphi}\ominus zM_{\theta,\varphi}$ and $M_{\theta,\varphi}\ominus wM_{\theta,\varphi}$, respectively. So by Proposition \ref{prop4.5},
\[ \Sigma_0^{M_{\theta,\varphi}}(0,0)=\sum_{m,n,k,l\geq 0} |\langle \beta_n(w) e_k(\theta,\varphi),
 \alpha_m(z) f_l(\theta,\varphi)\rangle|^2. \]

Assume $B(0,0)=(a,b)$, and write
\[1=\sum_{m,i\geq 0}\langle 1, \alpha_m  \theta^i \rangle \cdot \alpha_m(z) \theta^i(z)=\sum_{m,i\geq 0} \overline{\alpha_m(0)a^i} \cdot \alpha_m(z)\theta^i(z) \]
and
\[1=\sum_{n,j\geq 0}\langle 1, \beta_n \varphi^j \rangle \cdot \beta_n(w) \varphi^j(w)=\sum_{n,j\geq 0} \overline{\beta_n(0)b^j} \cdot \beta_n(w)\varphi^j(w). \]
By a direct computation, $\Sigma_0^{M_{\theta,\varphi}}(0,0)=$
\begin{align*}
&\sum_{m,n,k,l\geq 0} \mid\langle \sum_{m',i\geq 0} \overline{\alpha_{m'}(0)a^i} \cdot \alpha_{m'}\beta_n \theta^i e_k(\theta,\varphi),
 \sum_{n',j\geq 0} \overline{\beta_{n'}(0)b^j} \cdot \alpha_m \beta_{n'}\varphi^j f_l(\theta,\varphi)\rangle \mid^2\\
 &=\sum_{m,n,k,l\geq 0} \mid \sum_{i,j\geq 0} \overline{\alpha_m(0)}\beta_n(0)   \langle \bar{a}^i z^i e_k(z,w), \bar{b}^j w^j f_l(z,w) \rangle \mid^2\\
 &=\sum_{m,n,k,l\geq 0} \mid  \overline{\alpha_m(0)}\beta_n(0)   \langle \frac{1}{1-\bar{a}z} e_k(z,w), \frac{1}{1-\bar{b}w} f_l(z,w) \rangle \mid^2\\
 &=\sum_{k,l\geq 0} \sum_{m,n\geq 0} |\alpha_m(0)|^2|\beta_n(0)|^2 \mid  \langle \frac{1}{1-\bar{a}z} e_k(z,w), \frac{1}{1-\bar{b}w} f_l(z,w) \rangle \mid^2\\
 &=(1-|\theta(0)|^2)(1-|\varphi(0)|^2)\sum_{k,l\geq 0} \mid  \langle \frac{1}{1-\bar{a}z} e_k(z,w), \frac{1}{1-\bar{b}w} f_l(z,w) \rangle \mid^2\\
 &=\Sigma_0^M(\theta(0),\varphi(0)).
\end{align*}
Similarly, again by Proposition \ref{prop4.5},
\begin{align*}
  \Sigma_1^{M_{\theta,\varphi}}(0,0)&=\sum_{m,n,k,l\geq 0} |\langle w\beta_n(w) e_k(\theta,\varphi),
 z\alpha_m(z) f_l(\theta,\varphi)\rangle|^2 \\
   & =\sum_{m,n,k,l\geq 0} |\langle \bar{z} \beta_n(w) e_k(\theta,\varphi),
 \bar{w} \alpha_m(z) f_l(\theta,\varphi)\rangle_{L^2(\mathbb{T}^2)}|^2.
\end{align*}
Write
\begin{align*}
L^2(\mathbb{T})&=\overline{H^2(\mathbb{T})} \bigoplus  zH^2(\mathbb{T})\\
               &=\left(\bigoplus_{i\geq 0} \overline{\theta^i \mathcal{K}_\theta}\right)\bigoplus \left(\bigoplus_{i\geq 0} z \theta^i \mathcal{K}_\theta \right)\\
               &=\left(\bigoplus_{i<0} \theta^i \overline{\mathcal{K}_\theta}\right)\bigoplus \overline{\mathcal{K}_\theta}\bigoplus z\mathcal{K}_\theta \bigoplus \left(\bigoplus_{i>0} z\theta^i \mathcal{K}_\theta \right),
\end{align*}
and
\begin{align*}
  \bar{z}&=\sum_{m\geq 0}\sum_{i=-\infty}^\infty \langle \bar{z}, \alpha_m  \theta^i \rangle_{L^2(\mathbb{T}^2)} \cdot \alpha_m(z) \theta^i(z)\\
   &=\sum_{m\geq 0}\sum_{i<0} \langle \theta^{-i}\overline{\alpha_m }, z\rangle  \alpha_m(z) \theta^i(z)\\
   &=\sum_{m\geq 0}\sum_{i\geq 0} \langle \theta^i (\theta \overline{\alpha_m }), z \rangle \cdot
    \alpha_m(z)\overline{\theta^{i+1}}.
\end{align*}
Since the multiplication operator $M_\theta$ on $L^2(\mathbb{T})$ is unitarily equivalent to a multiple bilateral shift, with its multiplicity determined by the dimension of model space $\mathcal{K}_\theta$, we get that $M_\theta$ is unitary and $\theta \overline{\mathcal{K}_\theta}=z\mathcal{K}_\theta$. Thus we can assume $\theta \overline{\alpha_m}=z\alpha'_m$ with $\{\alpha'_m\}_{m\geq 0}$  remains an orthonormal basis of $\mathcal{K}_\theta$. In this case, we have
\[ \bar{z}=\sum_{m\geq 0}\sum_{i\geq 0} \langle \theta^i \alpha'_m, 1 \rangle \cdot
    \alpha_m(z)\overline{\theta^{i+1}} =\sum_{m, i\geq 0} a^i \alpha'_m(0) \cdot \alpha_m(z)\overline{\theta^{i+1}} .  \]
Similarly, assume $\varphi \overline{\beta_n}= w\beta'_n$ with $\{\beta'_n\}_{n\geq 0}$  remains an orthonormal basis of $\mathcal{K}_\varphi$. Then we can write
\[ \bar{w}=\sum_{n\geq 0}\sum_{j\geq 0} \langle \varphi^j \beta'_n, 1 \rangle \cdot
    \beta_n(z)\overline{\varphi^{j+1}} =\sum_{n, j\geq 0} b^j \beta'_n(0) \cdot \beta_n(z)\overline{\varphi^{j+1}}. \]
Thus we have $\Sigma_1^{M_{\theta,\varphi}}(0,0)=$
\begin{align*}
&\sum_{m,n,k,l\geq 0} \mid\langle \sum_{m', i \geq 0 }  a^{i}\alpha'_{m'}(0)  \cdot \alpha_{m'}\beta_n  \overline{\theta^{i+1}} e_k(\theta,\varphi),
 \sum_{n', j\geq 0} b^j \beta'_{n'}(0)  \cdot \alpha_m \beta_{n'}\overline{\varphi^{j+1}} f_l(\theta,\varphi)\rangle_{L^2(\mathbb{T}^2)} \mid^2\\
 &=\sum_{m,n,k,l\geq 0} \mid \sum_{i,j\geq 0} \langle  a^i \alpha'_m(0) \alpha_m \beta_n \varphi^{j+1} e_k(\theta,\varphi), b^j \beta'_n(0) \alpha_m \beta_n \theta^{i+1}f_l(\theta,\varphi)     \rangle  \mid^2\\
 &=\sum_{m,n,k,l\geq 0} \mid \sum_{i,j\geq 0} \langle  a^i \alpha'_m(0) w^{j+1} e_k(z,w), b^j \beta'_n(0) z^{i+1} f_l(z,w)     \rangle  \mid^2\\
  &=\sum_{m,n,k,l\geq 0} \mid  \alpha'_m(0)\overline{\beta'_n(0)}   \langle \frac{w}{1-\bar{b}w} e_k(z,w), \frac{z}{1-\bar{a}z} f_l(z,w) \rangle \mid^2\\
 &=\sum_{k,l\geq 0} \sum_{m,n\geq 0} |\alpha'_m(0)|^2|\beta'_n(0)|^2 \mid  \langle \frac{w}{1-\bar{b}w} e_k(z,w), \frac{z}{1-\bar{a}z} f_l(z,w) \rangle \mid^2\\
 &=(1-|\theta(0)|^2)(1-|\varphi(0)|^2)\sum_{k,l\geq 0} \mid  \langle \frac{w}{1-\bar{b}w} e_k(z,w), \frac{z}{1-\bar{a}z} f_l(z,w) \rangle \mid^2\\
 &=\Sigma_1^M(\theta(0),\varphi(0)).
\end{align*}
\end{proof}

\begin{corollary}\label{thm4.8}
Let  $\theta,\varphi$ be two nonconstant inner functions and $M$ be a Hilbert-Schmidt submodule, then
\[ \Sigma_0^{M_{\theta,\varphi}}(a,b)=\Sigma_0^M(\theta(a),\varphi(b)),~~~ \Sigma_1^{M_{\theta,\varphi}}(a,b)=\Sigma_1^M(\theta(a),\varphi(b))\]
\end{corollary}
\begin{proof}
By Theorem \ref{thm4.6} and the fact that $C_{\varphi_a,\varphi_b}M_{\theta,\varphi}=M_{\theta\circ \varphi_a,~\varphi\circ\varphi_b}$,
\[ \Sigma^{M_{\theta,\varphi}}_0(a,b)=\Sigma_0^{C_{\varphi_a,\varphi_b}M_{\theta,\varphi}}(0,0)=
\Sigma_0^{M_{\theta\circ \varphi_a,~\varphi\circ\varphi_b}}(0,0)=\Sigma_0^M(\theta(a),\varphi(b)).
 \]
The other formula follows from a similar discussion.
\end{proof}

Combine with the fact that $\Sigma^{M}_0(0,0)-\Sigma^{M}_1(0,0)=1$ for any Hilbert-Schmidt submodule, we get
\begin{corollary}\label{ind}
If $M$ is a Hilbert-Schmidt submodule, then
\[ \Sigma^{M}_0(a,b)-\Sigma^{M}_1(a,b)\equiv 1. \]
Thus for any $(a,b)\in \mathbb{D}^2$, $\ind F_{a,b}=-1$.
\end{corollary}
\begin{proof}
\[ \Sigma^{M}_0(a,b)-\Sigma^{M}_1(a,b)=\Sigma_0^{M_{\varphi_a,\varphi_b}}(0,0)-\Sigma_1^{M_{\varphi_a,\varphi_b}}(0,0)=1.  \]
\end{proof}
\begin{remark}
If $M$ is a  Hilbert-Schmidt submodule such that the core operator $C_M=P_0Q_0-P'_0Q'_0$ is trace class, then we can prove that
\[\Sigma^M_0(a,b )- \Sigma^M_1(a,b )\equiv 1\]
in another simple way.
\begin{proof}
Recall that the reproducing kernel for $H^2(\mathbb{D}^2)$ is the Szeg\"{o}  kernel $K_{\lambda,\mu}(z,w)=\frac{1}{1-\bar{\lambda}z} \frac{1}{1-\bar{\mu}w}$. For every bounded operator $A$ on $H^2(\mathbb{D}^2)$, we define the \textbf{Berezin transform} of $A$ as the function
\[ A(\lambda,\mu)=\langle A K_{\lambda,\mu},K_{\lambda,\mu}\rangle,~~~~~(\lambda,\mu)\in \mathbb{D}^2. \]
The function $A(\lambda,\mu)$ is said to have $L^1$ boundary value if $A(rz,rw)$ converges in the $L^1(\mathbb{T}^2)$ norm when $r$ increases to $1$. For simplicity, the boundary value function $A(z,w)$ is also denoted by $A(\lambda,\mu)$ in contexts where no confusion is likely to arise. 
Then the following classical result is essential for our proof.
\begin{theorem}[Theorem 2.4,  \cite{Ya2004}]\label{trace}
If $A$ is trace class, then $A(\lambda,\mu)$ has $L^1$ boundary value and
\[ \tr(A)=\int_{\mathbb{T}^2} A(\lambda,\mu)|d\lambda||d\mu|. \]
\end{theorem}
It is easy to check that
\begin{align*}
  P_aK_{\lambda,\mu}^M(z,w) &=(1-\overline{\varphi_a(\lambda)}\varphi_a(z))K_{\lambda,\mu}^M(z,w), \\
  Q_bK_{\lambda,\mu}^M(z,w) &=(1-\overline{\varphi_b(\mu)}\varphi_b(w))K_{\lambda,\mu}^M(z,w), \\
P'_aK_{\lambda,\mu}^M(z,w) &=
\overline{\varphi_a(\lambda)}\varphi_a(z)(1-\overline{\varphi_b(\mu)}\varphi_b(w))K_{\lambda,\mu}^M(z,w), \\
Q'_bK_{\lambda,\mu}^M(z,w) &=
\overline{\varphi_b(\mu)}\varphi_b(w)(1-\overline{\varphi_a(\lambda)}\varphi_a(z))K_{\lambda,\mu}^M(z,w).
\end{align*}
By the reproducing property of $K^M$, we have
\begin{align*}
 \langle P_aK_{\lambda,\mu}^M(z,w), Q_bK_{\lambda,\mu}^M(z,w)\rangle & =(1-|\varphi_a(\lambda)|^2-|\varphi_b(\mu)|^2) K^M_{\lambda,\mu}(\lambda,\mu)\\
 &+\langle \overline{\varphi_a(\lambda)}\varphi_a(z)K^M_{\lambda,\mu}(z,w), \overline{\varphi_b(\mu)}\varphi_b(w)K^M_{\lambda,\mu}(z,w) \rangle,\\
 \langle P'_aK_{\lambda,\mu}^M(z,w), Q'_bK_{\lambda,\mu}^M(z,w)\rangle & =\langle \overline{\varphi_a(\lambda)}\varphi_a(z)K^M_{\lambda,\mu}(z,w), \overline{\varphi_b(\mu)}\varphi_b(w)K^M_{\lambda,\mu}(z,w) \rangle\\
 &-|\varphi_a(\lambda)|^2|\varphi_b(\mu)|^2 K^M_{\lambda,\mu}(\lambda,\mu).
\end{align*}
Thus by Theorem \ref{trace},
\begin{align*}
 \Sigma^M_0(a,b )- \Sigma^M_1(a,b ) &=\lim_{r\rightarrow 1^-} \int_{\mathbb{T}^2} \langle P_aK_{r\lambda,r\mu}^M, Q_bK_{r\lambda,r\mu}^M\rangle - \langle P'_aK_{r\lambda,r\mu}^M, Q'_bK_{r\lambda,r\mu}^M\rangle |d\lambda| |d\mu| \\
  &=\lim_{r\rightarrow 1^-} \int_{\mathbb{T}^2} (1-|\varphi_a(r\lambda)|^2)(1-|\varphi_b(r\mu)|^2)K^M(r\lambda,r\mu)(r\lambda,r\mu)  |d\lambda| |d\mu|\\
  &=\lim_{r\rightarrow 1^-} \int_{\mathbb{T}^2} \frac{(1-|\varphi_a(r\lambda)|^2)(1-|\varphi_b(r\mu)|^2)}{(1-|r\lambda|^2)(1-|r\mu|^2)} G^M_{r\lambda,r\mu}(r\lambda,r\mu)   |d\lambda| |d\mu|.
\end{align*}
Since the core function $G^M_{\lambda,\mu}(\lambda,\mu)$ has boundary value $1$ at almost every point in $\mathbb{T}^2$, then
\begin{align*}
  \Sigma^M_0(a,b )- \Sigma^M_1(a,b ) &=\int_{\mathbb{T}^2} \frac{(1-|\varphi_a(\lambda)|^2)(1-|\varphi_b(\mu)|^2)}{(1-|\lambda|^2)(1-|\mu|^2)} |d\lambda| |d\mu| \\
   & =\int_{\mathbb{T}^2} \frac{1-|a|^2}{|1-\bar{a}\lambda|^2}\frac{1-|b|^2}{|1-\bar{b}\mu|^2}  |d\lambda| |d\mu|\\
   &=1.
\end{align*}
\end{proof}
\end{remark}

\begin{remark}

Once we know the fact that $\Sigma^M_0(a,b )- \Sigma^M_1(a,b )\equiv 1$, then Theorem \ref{thm4.6}  follows readily from the equation $G^{M_{\theta,\varphi}}_{\lambda,\mu}(z,w)=G^{M}_{\theta(\lambda),\varphi(\mu)}(\theta(z),\varphi(w))$.

\begin{proof}
Let  $\theta,\varphi$ be two nonconstant one-variable inner functions and set $(a,b)=(\theta(0),\varphi(0))$, we first show that
\begin{equation}\label{4.4}
  \|G^{M_{\theta,\varphi}}\|_{L^2(\mathbb{T}^2\times\mathbb{T}^2)}^2=
\|G^{M_{\varphi_a,\varphi_b}}\|_{L^2(\mathbb{T}^2\times\mathbb{T}^2)}^2 .
\end{equation}
In fact, by Theorem \ref{Shapiro}, given an analytic or co-analytic function  $f$ on the disk and an analytic self mapping $\varphi$ of the disk, the $2$-norm of $f\circ\varphi$ is only depended on $f$ and the Nevanlinna counting function $N_\varphi(w)$. Also, by Lemma \ref{Nevanlina}, we know that when $\varphi$ is inner function, the Nevanlinna counting function $N_\varphi(w)$ is determined by $\varphi(0)$. Hence for any analytic or co-analytic function  $f$, if $\varphi$ and $\psi$ are two inner functions with $\varphi(0)=\psi(0)$, then
$\|f\circ\varphi\|_2^2=\|f\circ\psi\|_2^2$. Since the core function $G^{M}_{\lambda,\mu}(z,w)$ is either co-analytic or analytic in each variable, then the conclusion follows directly from Corollary \ref{Cor3.1}.

By a direct computation or Proposition 1.3 in \cite{Ya2005}, one can check that
\[ C_{M_{\varphi_a,\varphi_b}}=C_{\varphi_a,\varphi_b}C_MC^*_{\varphi_a,\varphi_b}=C_M(a,b).\] 
Thus
\[\|G^{M_{\varphi_a,\varphi_b}}\|_{L^2(\mathbb{T}^2\times\mathbb{T}^2)}^2=\|C_{M_{\varphi_a,\varphi_b}}\|_{H.S.}^2
=\Sigma_0^{M}(a,b)+\Sigma_1^{M}(a,b)=2\Sigma_0^{M}(a,b)-1.\]
Then Theorem \ref{thm4.6} follows directly from  formula (\ref{4.4}) and the fact that
\[\|G^{M_{\theta,\varphi}}\|_{L^2(\mathbb{T}^2\times\mathbb{T}^2)}^2=\Sigma_0^{M_{\theta,\varphi}}(0,0)
+\Sigma_1^{M_{\theta,\varphi}}(0,0)=2\Sigma_0^{M_{\theta,\varphi}}(0,0)-1.\]
\end{proof}
\end{remark}

A more tempting problem concerns the connection between the submodule $M$ and the numerical invariant functions.
\begin{problem}
Let $M$ be an arbitrary Hilbert-Schmidt submodule, what are the common properties of numerical invariant function $\Sigma^M_0(a,b)$ or $\Sigma^M_1(a,b)$? For example, is $\Sigma^M_0(a,b)$ bounded for all Hilbert-Schmidt submodules? Moreover, what information of submodule $M$ can we learn from the numerical invariant functions?
\end{problem}

\section{The Submodule $[z-w]$}
In \cite{z-w}, K. Guo, S. Sun, D. Zheng and C. Zhong developed a machinery to study multiplication operators on the Bergman space via the quotient module $[z-w]^\perp$. Moreover, by using this machinery, S. Sun and D. Zheng \cite{Sun} gave a new proof of the Beurling-type theorem for the Bergman space. This machinery has been shown to be very useful in dealing with some question on Bergman space so far. In this section, we will show that $\Sigma^{[z-w]}_1(a,b)$ is bounded by $2$. For convenience, we omit the superscript of $\Sigma^{M}_i(a,b)$.

Set
\[e_n(z,w)=\frac{1}{\sqrt{n+1}}\cdot \frac{z^{n+1}-w^{n+1}}{z-w}=\frac{1}{\sqrt{n+1}}\sum_{i=0}^{n}z^{n-i}w^i,~~~n=0,1,2,\cdots.\]
One verifies that $\{e_n\mid n\geq 0\}$ is an orthonormal basis for quotient module$[z-w]^\perp$.

For submodule $M=[z-w]$, we define
\[ \phi_n(z,w)=\frac{1}{\sqrt{n+2}}(ze_n(z,w)-\sqrt{n+1}w^{n+1}),~~~ n\geq 0 \]
and
\[ \psi_n(z,w)=\phi_n(w,z)=\frac{1}{\sqrt{n+2}}(we_n(z,w)-\sqrt{n+1}z^{n+1}), ~~~n\geq 0.  \]
With some computation, one can verify that $\{\phi_n\mid n\geq 0\}$ is an orthonormal basis for $M\ominus zM$ and  $\{\psi_n\mid n\geq 0\}$ is an orthonormal basis for $M\ominus wM$ (see Example 3(c) in \cite{Ya3}).

\begin{theorem}{\label{thm5.3}}
For the submodule $M=[z-w]$ and any $(a,b)\in \mathbb{D}^2$, we have
\[ \Sigma_1(a,b)\leq 2. \]
\end{theorem}
\begin{proof}
First note that
\begin{equation}\label{6.2}
  \langle w^{n+1}\phi_k, z^{m+1}\psi_l \rangle =\frac{1}{\sqrt{k+2}}\cdot \frac{1}{\sqrt{l+2}}\langle w^n e_k , z^m e_l  \rangle,
\end{equation}
then by the formula \eqref{4.6} in Proposition \ref{prop4.5},
\begin{align*}
  \Sigma_1(a,b)&=(1-|a|^2)(1-|b|^2)\sum_{k,l\geq 0} |\langle \frac{w \phi_k}{1-\bar{b}w}, \frac{z\psi_l}{1-\bar{a}z} \rangle|^2 \\
   & =(1-|a|^2)(1-|b|^2)\sum_{k,l\geq 0} \frac{1}{(k+2)(l+2)} \cdot |\sum_{m,n\geq 0} a^m \bar{b}^n  \langle w^n e_k , z^m e_l \rangle |^2.
\end{align*}
By some computations, we have
\begin{equation*}
|\sum_{m,n\geq 0} a^m \bar{b}^n \langle w^n e_k , z^m e_l \rangle |^2=
\begin{cases}
  \frac{|a|^{2(l-k)}}{(k+1)(l+1)}~|\sum_{i=0}^k (k+1-i) (a\bar{b})^i|^2&if\quad k<l\\
  \frac{1}{(k+1)^2}~|\sum_{i=0}^k (k+1-i) (a\bar{b})^i|^2&if\quad k=l\\
  \frac{|b|^{2(k-l)}}{(k+1)(l+1)}~|\sum_{i=0}^l (l+1-i) (a\bar{b})^i|^2&if\quad k>l\\
\end{cases}
\end{equation*}
Thus it is easy to see that $ \Sigma_1(a,b)\leq  \Sigma_1(|a|,|b|)$.

Without loss of generality, we assume $0\leq a,b < 1$, and we have
\begin{align*}
  \Sigma_1(a,b) &=\sum_{k\geq 0} \frac{(1-a^2)(1-b^2)}{(k+1)^2(k+2)^2}|\sum_{i=0}^k (k+1-i) (ab)^i|^2 \\
   &+\sum_{l>k\geq 0}\frac{(1-a^2)(1-b^2)}{(l+1)(l+2)}\cdot \frac{a^{2(l-k)}}{(k+1)(k+2)} |\sum_{i=0}^k (k+1-i) (ab)^i|^2\\
   &+\sum_{l>k\geq 0}\frac{(1-a^2)(1-b^2)}{(l+1)(l+2)}\cdot \frac{ b^{2(l-k)}}{(k+1)(k+2)}|\sum_{i=0}^k (k+1-i) (ab)^i|^2\\
   &\leq \sum_{l\geq k\geq 0}\frac{(1-a^2)(1-b^2)}{(l+1)(l+2)}\cdot \frac{a^{2(l-k)}+b^{2(l-k)}}{(k+1)(k+2)} |\sum_{i=0}^k (k+1-i) (ab)^i|^2\\
   &=\sum_{m, k\geq 0}\frac{a^{2m}+b^{2m}}{(m+k+1)(m+k+2)}\cdot \frac{(1-a^2)(1-b^2)}{(k+1)(k+2)} |\sum_{i=0}^k (k+1-i) (ab)^i|^2\\
\end{align*}
It is not hard to check that
\begin{equation}\label{6.3}
  \sum_{i=0}^k (k+1-i) (ab)^i= \sum_{i=1}^{k+1} \frac{1-(ab)^i}{1-ab}.
\end{equation}

\begin{lemma}\label{lemma6.4}
For $0\leq a,b < 1$ and $i\geq 1$.
$$ (1-a^ib^i)^2(1-a^2)(1-b^2)\leq (1-ab)^2(1-a^{2i})(1-b^{2i}).$$
\end{lemma}
\begin{proof}
\begin{align*}
  \left( \frac{1-a^ib^i}{1-ab}\right)^2 &=\left( \sum_{k=0}^{i-1} a^k b^k \right)^2 \leq \left( \sum_{k=0}^{i-1} a^{2k} \right)\left( \sum_{k=0}^{i-1}  b^{2k} \right)=\frac{1-a^{2i}}{1-a^2}\cdot \frac{1-b^{2i}}{1-b^2}.
\end{align*}
\end{proof}
By Lemma \ref{lemma6.4},
\[ \sqrt{(1-a^2)(1-b^2)}\cdot \frac{1-a^ib^i}{1-ab} \leq \sqrt{(1-a^{2i})(1-b^{2i})}. \]
Combine with formula \eqref{6.3}, we have
\begin{align*}
  (1-a^2)(1-b^2) | \sum_{i=0}^k (k+1-i) (ab)^i |^2 &=|\sum_{i=1}^{k+1} \sqrt{(1-a^2)(1-b^2)}\cdot \frac{1-a^ib^i}{1-ab} |^2 \\
   & \leq |\sum_{i=1}^{k+1}  \sqrt{(1-a^{2i})(1-b^{2i})} |^2\\
   &\leq \left(\sum_{i=1}^{k+1} 1-a^{2i} \right)\left(\sum_{i=1}^{k+1} 1-b^{2i} \right)\\
   &\leq (k+1)^2 (1-a^{2(k+1)})(1-b^{2(k+1)}).
\end{align*}
Thus
\begin{align*}
  \Sigma_1(a,b) & = \sum_{m, k\geq 0}\frac{a^{2m}+b^{2m}}{(m+k+1)(m+k+2)}\cdot \frac{(1-a^2)(1-b^2)}{(k+1)(k+2)} |\sum_{i=0}^k (k+1-i) (ab)^i|^2 \\
   & \leq \sum_{m, k\geq 0} \frac{a^{2m}+b^{2m}}{(m+k+1)(m+k+2)}\cdot \frac{(k+1)^2}{(k+1)(k+2)} (1-a^{2(k+1)})(1-b^{2(k+1)})\\
  &\leq \sum_{m\geq 0, k\geq 1}\frac{(a^{2m}+b^{2m})(1-a^{2k})(1-b^{2k})}{(m+k)(m+k+1)}.
\end{align*}

\begin{lemma}\label{lemma6.5}
For $0\leq a< 1$,
$$\sum_{ m\geq 0, k\geq 1} \frac{(1-a^k)a^m}{(m+k)(m+k+1)}=1 $$
\end{lemma}
\begin{proof}
\begin{align*}
  \sum_{k\geq 1, m\geq 0} \frac{(1-a^k)a^m}{(m+k)(m+k+1)} & =\sum_{k\geq 1, m\geq 0} \frac{a^m}{(m+k)(m+k+1)}-\sum_{k\geq 1, m\geq 0} \frac{a^{k+m}}{(m+k)(m+k+1)} \\
   & =\sum_{m\geq 0} \frac{a^m}{m+1}-\sum_{n(=m+k)\geq 1} \sum_{k=1}^n \frac{a^{n}}{n(n+1)}\\
   &=\sum_{m\geq 0} \frac{a^m}{m+1}-\sum_{n\geq 1}\frac{a^{n}}{n+1}\\
   &=1.
\end{align*}
\end{proof}

Therefore
\begin{align*}
  \Sigma_1(a,b) & \leq \sum_{m\geq 0, k\geq 1}\frac{(a^{2m}+b^{2m})(1-a^{2k})(1-b^{2k})}{(m+k)(m+k+1)}.  \\
   & \leq \sum_{m\geq 0, k\geq 1}\frac{a^{2m}(1-a^{2k})+b^{2m}(1-b^{2k})}{(m+k)(m+k+1)}\\
   &=2
\end{align*}
\end{proof}

As a direct consequence of Theorem \ref{thm4.6}, Corollary \ref{ind} and Theorem \ref{thm5.3}, we obtain the following result.
\begin{corollary}
Let  $\theta(z),\varphi(w)$ be two nonconstant inner functions, then the submodule $[\theta(z)-\varphi(w)]$ is Hilbert-Schmidt and
\[  \|C_{[\theta(z)-\varphi(w)]}\|_{H.S.}^2=2\Sigma_1(\theta(0),\varphi(0))+1 \leq 5 .  \]
\end{corollary}
It is not hard to see that our estimate is not optimal, and it's natural to ask the following question:
\begin{question}
For the submodule $[z-w]$, does it holds that
\[ \Sigma_1(a,b)\leq \Sigma_1(0,0)=\frac{\pi^2}{6}-1?\]
\end{question}

\vspace{0.5cm}
{\bf Acknowledgements} The authors sincerely thank the anonymous reviewers for their helpful feedback and constructive suggestions, which greatly improved this paper. The authors also appreciate Prof. Rongwei Yang for many valuable discussions and insightful feedbacks, which greatly contributed to the development of this work. The first author is supported by National Natural Science Foundation of China (No.12401151), and the Postdoctoral Researcher Foundation of China (Grant No.GZB20240100). The second author is supported by National Natural Science Foundation of China (No.12031002).

\vspace{0.5cm}
{\bf Data Availability Statement} No data was used for the research described in the article.

\end{document}